\documentclass{ifacconf}

\usepackage{amssymb}
\usepackage{amsmath}

\usepackage{graphicx}      
\usepackage{natbib}        
\usepackage{algorithm}
\usepackage{algpseudocode}

\newcommand{\norm}[1]{\left\lVert#1\right\rVert}

\begin{document}
\begin{frontmatter}

\title{Reinforcement Learning based Design of Linear Fixed Structure Controllers}

\thanks[footnoteinfo]{\textcopyright 2020 the authors. This work has been accepted to IFAC World Congress for publication under a Creative Commons Licence CC-BY-NC-ND}

\author[First]{Nathan P. Lawrence} 
\author[Third]{Gregory E. Stewart} 
\author[First]{Philip D. Loewen}
\author[Fourth]{Michael G. Forbes}
\author[Fourth]{Johan U. Backstrom}
\author[Second]{R. Bhushan Gopaluni}
\address[First]{Department of Mathematics, University of British Columbia, Vancouver, BC V6T 1Z2, Canada (e-mail: lawrence@math.ubc.ca, loew@math.ubc.ca).}
\address[Second]{Department of Chemical and Biological Engineering, University of British Columbia, Vancouver, BC V6T 1Z3, Canada (e-mail: bhushan.gopaluni@ubc.ca)}
\address[Third]{Department of Electrical and Computer Engineering, University of British Columbia, Vancouver, BC V6T 1Z4, Canada (e-mail: stewartg@ece.ubc.ca)}
\address[Fourth]{Honeywell Process Solutions, North Vancouver, BC V7J 3S4, Canada (e-mail: michael.forbes@honeywell.com, johan.backstrom@honeywell.com)}

\begin{abstract}                
Reinforcement learning has been successfully applied to the problem of tuning PID controllers in several applications. The existing methods often utilize function approximation, such as neural networks, to update the controller parameters at each time-step of the underlying process. In this work, we present a simple finite-difference approach, based on random search, to tuning linear fixed-structure controllers. For clarity and simplicity, we focus on PID controllers. Our algorithm operates on the entire closed-loop step response of the system and iteratively improves the PID gains towards a desired closed-loop response. This allows for embedding stability requirements into the reward function without any modeling procedures. 
\end{abstract}


\end{frontmatter}

\section{Introduction}

Reinforcement learning (RL) is a branch of machine learning in which the objective is to learn an optimal strategy for interacting with an environment through experiences \citep{sutton2018reinforcement}. Traditional tabular methods for RL do not apply in continuous state or action spaces and are cumbersome in high-dimensional settings. For instance, the game of Go contains an intractable number of possible board configurations, which motivates the synthesis of deep learning with RL \citep{silver2016mastering}. 

The success of RL methods reported in the literature is due to increasingly complicated algorithms. Combined with the inherent stochasticity due to random seeds or the underlying environment itself, as well as sensitivity to hyperparameters, the problem of reproducibility has become prevalent \citep{henderson2018deep, islam2017reproducibility}. Several recent works have proposed simple algorithms that achieve performance competitive or superior to standard $Q$-learning and policy gradient methods \citep{salimans2017evolution, rajeswaran2017towards, mania2018simple}.

The applications of machine learning and RL to process control are relatively recent and limited in industrial implementation \citep{venkatasubramanian2019promise, spielberg2019toward}. Among the first approaches were \cite{lee2001neuro, lee2008value}, in which the authors develop an approximate dynamic programming approach with function approximation as a computationally efficient framework for model predictive control and gain scheduling, respectively. More recently, \cite{spielberg2019toward} and \cite{wang2018novel} proposed deep RL algorithms for control of discrete-time nonlinear processes. Both approaches are in the class of actor-critic methods, in which the actor (controller) is represented by a deep neural network. It is then shown empirically that networks represent a flexible class of controllers capable of learning to control complex systems. In contrast, PID controllers are simple and industrially established mechanisms for set-point tracking. However, setting the PID tuning parameters is known to be challenging and represents a nonlinear design problem. PID tuning thus represents both a challenging RL problem, and one which has the practical goal of being implemented in a production control system.

\begin{figure}[thb]
\begin{center}
\includegraphics[width=8.4cm]{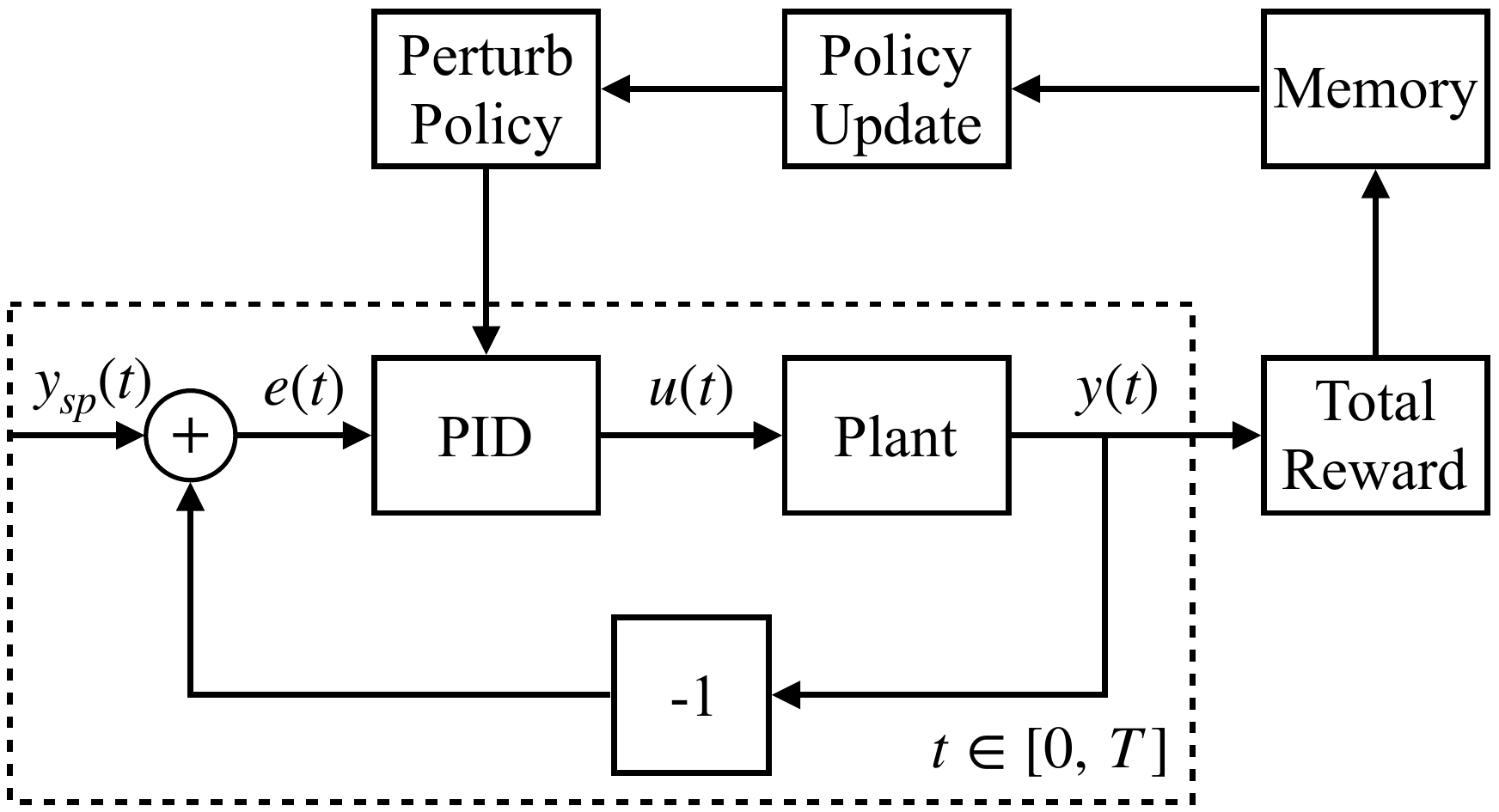}    
\caption{A standard closed-loop structure is shown inside the dashed box. Arrows crossing the dashed line indicate the passing of some time-horizon $[0,\ T]$. Outside the dashed box, we store cumulative rewards based on slightly perturbed policies, which are used to update the policy with the finite-difference scheme described in section \ref{sec:algorithm}.} 
\label{fig:AlgSchematic}
\end{center}
\end{figure}

In this work, we develop an adaptive tuning algorithm based on a simple random search procedure for linear fixed-structure controllers. For simplicity in the development we focus on PID controllers. Our algorithm tracks a desired closed-loop step response by evaluating the distance between the desired response and the response generated by slightly perturbing the policy which produces the controller parameters. We update the policy using a finite-difference approximation of the objective. Finally, although our method does not make use of a plant model, we focus on single-input single-output systems.

This paper is organized as follows: Section \ref{sec:PID} gives a brief description of PID control. Section \ref{sec:RL} outlines common methods in reinforcement learning. Section \ref{sec:PID_RL} describes our approach and algorithm for PID tuning based on simplified RL strategies discussed in section \ref{sec:RL}. Further, we contrast our approach to other RL-based tuning approaches. Finally, we show several simulation results in section \ref{sec:results}.

\section{PID control}\label{sec:PID}

In this section, we highlight some common strategies for PID control as they motivate our approach presented in the following section.

We use the parallel form of the PID controller:
\begin{equation}
    u(t) = k_p e(t) + k_i \int_{0}^{t} e(\tau) d \tau + k_d \frac{d}{d t} e(t).
\label{eq:PID}
\end{equation}

Although the structure of a PID controller is simple, requiring only three tuning parameters $k_p, k_i, k_d$, adjusting these parameters to meet certain performance specifications is difficult in practice. Below we describe some performance metrics and strategies for tuning $k_p, k_i, k_d$.

\subsection{Performance measures}

Our proposed algorithm in section \ref{sec:PID_RL} does not rely on a plant model. Therefore, to evaluate the performance of the closed-loop step response of a system we use the accumulated error over some time horizon $[0, T]$. Common measures include the integral error, such as integral absolute error (IAE) or integral squared error (ISE):
\begin{equation}
    \text{IAE } = \int_{0}^{T} \left|e(t)\right| d t \quad \text{and} \quad \text{ISE } = \int_{0}^{T} e(t)^2 d t.\label{eq:IAE}
\end{equation}

Note that in order for the IAE or ISE to be a useful measure of performance, $T$ should be large enough to measure the error accumulated through the closed-loop step response up to steady-state. In practice, we approximate the IAE and ISE through sampling at discrete time steps.

The criteria in \eqref{eq:IAE} motivate our reward function in the following sections. However, IAE and ISE may also incorporate a weight on the magnitude of the control signal or its first-order difference.

\subsection{Internal Model Control}\label{sec:IMC}

Internal model control (IMC) utilizes a reference model in the feedback loop by incorporating the deviation of the plant output from that of the model. The resulting control structure often results in a PID controller for a large number of single-input single-output processes \citep{rivera1986internal}. The simple PID tuning rules of \cite{skogestad2001probably} provide robust performance based on a first-order or second-order approximate model of the plant. Such a model can be obtained from a single step test of the plant. One then measures the plant gain, and time delay, as well as first and second order time constants. The resulting model can then be used to derive suitable PID gains for the plant using, for example, the SIMC rules \citep{skogestad2001probably}. In section \ref{sec:results}, we use SIMC to initialize the PID gains in algorithm \ref{alg:MainAlg}. 

\section{Reinforcement Learning}\label{sec:RL}

In this section, we provide a brief overview of common RL methods and their applications in process control. We then formulate the problem of tuning a PID controller as a RL problem and describe our approach to solving it. 

\subsection{The Reinforcement Learning Problem}\label{subsec:RLprob}

For each 
\emph{state}\footnote{The `state' $x_i$ is some measurement which characterizes the environment.}
$x_i$ the agent encounters, it takes some \emph{action} $a_i$, leading to a new state $x_{i+1}$. Upon taking action $a_i$, the agent receives a reward $r(x_i, a_i)$. The reward measures how desirable a particular state-action interaction is.
To interact optimally with the environment the agent learns to maximize the cumulative reward due to a sequence of interactions. Formally, the environment is modeled by a Markov decision process with initial distribution $p(x_1)$ and transition probability $p(x_{i+1} | x_i, a_i)$. The agent then transitions from states to actions based on a conditional probability distribution $\pi$ referred to as a \emph{policy}. If $h = (x_1, a_1, r_1, , \ldots, x_j, a_j, r_j)$ is some trajectory generated by the policy $\pi$ with sub-sequential states distributed according to $p$, we write $h \sim p^{\pi}(\cdot)$. If we assume the policy has some parametric structure given by $\pi_{\theta}$, then the problem we aim to solve is:
\begin{equation}
\begin{aligned}
    &\text{maximize} && \mathbb{E}_{h \sim p^{\pi_{\theta}}(\cdot)}\big[ R(h) | x_0 \big]\\
    &\text{over all} && \theta \in \mathbb{R}^{n},
\end{aligned}
\label{eq:objective}
\end{equation}
where $R$ denotes the accumulated reward received over an arbitrary trajectory.

Common approaches to solving \eqref{eq:objective} involve $Q$-learning and the policy gradient theorem (see \cite{spielberg2019toward} and the references therein). In particular, variations of these methods involving function approximation have made RL problems tractable in continuous state and action spaces (\cite{lillicrap2015continuous, silver2014deterministic, sutton2000policy}). Otherwise, discretization of the state and action spaces is necessary, leading to the ``curse of dimensionality". These methods and variations thereof have led to remarkable results in the game of Go and simulated continuous control such as locomotion tasks in MuJoCo \citep{silver2016mastering, lillicrap2015continuous, todorov2012mujoco}. Methods involving function approximation are better suited for these applications than discretization due to the huge number of possible state-action configurations.

\subsection{Simple Reinforcement Learning Strategies}

In contrast to the significant results described above, the issue of \emph{reproducibility} has recently been examined in the RL community. \cite{islam2017reproducibility} and \cite{henderson2018deep} compare performance of open-source implementations of popular policy gradient algorithms on various MuJoCo tasks. Factors such as hyperparameters, neural network architecture, or random seeds can have a drastic effect on algorithm performance. 

In light of these issues, several recent works have proposed simpler algorithms with competitive or superior performance in benchmark MuJoCo tasks compared to the state-of-the-art performance reported with policy gradient methods.

\cite{salimans2017evolution} consider a gradient-free optimization method called \emph{evolution strategies} (ES) as the basis for solving the optimization problem \eqref{eq:objective}. The underlying idea of ES is to perturb the parameters of the policy, evaluate the policy, then combine the policy parameters with the best return. ES is simpler than the approaches highlighted in section \ref{subsec:RLprob} because there is no value function approximation procedure and its policy updates do not rely on computing the gradient of \eqref{eq:objective}.

While \cite{salimans2017evolution} show that ES is competitive with standard RL benchmarks, their approach uses a neural network parameterization for the policy and includes several back-end algorithmic enhancements. To this end, \cite{rajeswaran2017towards} achieve benchmark performance on standard MuJoCo locomotion tasks using a natural policy gradient algorithm with linear policies, thereby showing that neural networks are not necessary for these tasks. A synthesis of these approaches is proposed by \cite{mania2018simple}, in which a gradient-free algorithm for training linear policies is shown to achieve roughly equal overall performance on these MuJoCo locomotion tasks. The proposed algorithm of \cite{mania2018simple} is then the basis of our approach.

\section{PID fine-tuning via Reinforcement Learning}\label{sec:PID_RL}

In this section, we outline our strategy for PID tuning via reinforcement learning and contrast it with previous such approaches.

\subsection{States, Actions, and Rewards}

In our approach, we define the actions of the RL agent to be a vector of PID gains and the state to be a discretization of the closed-loop step response under these PID gains over a finite time horizon. Similarly, the target transfer function is represented by a vector of \emph{target data}---a compatible discretization of the target transfer function step response. In principle, the target data may come from a simulated plant even if the RL algorithm is operating on a physical system. The reward for a state-action pair is then the mean absolute (or squared) error between the state and the target data. Concretely, we write $y(t)$ for the value of the closed-loop step response at the particular time $t$. Then we choose a
sample count $n$ and a vector of sampling times $T = [0, t_1, \ldots,t_{n-1}]$, 
and write $x = [y(0), \ldots, y(t_{n-1})]$.
The corresponding target data $\tilde{y}(t_i)$ is contained in a vector $\tilde{x}$.
Finally, the action is denoted by $K = [k_p, k_i, k_d]$. 
We express the reward for action $K$ in state $x$ as 
\begin{equation}
    r(x, K) = - \frac{1}{n} \norm{x - \tilde{x}}_{q}^{q}
    \label{eq:reward}
\end{equation}
where the exponent $q \in \{1,2\}$ is fixed.

\subsection{Closed-Loop Transfer Function Tracking Algorithm}\label{sec:algorithm}

For Algorithm \ref{alg:MainAlg} we use the notation above. 
We refer to the process being controlled as $P$,
and introduce functions
$\texttt{feedback}(P, K)$ to describe the negative feedback loop for plant $P$ under the PID controller with given gains $K$, 
and $\texttt{step}(H , T)$ to generate the vector of outputs at each time-step in $T$ for some given system $H$.
The goal is to determine a $k\times n$ matrix $M$ for which optimal PID gains $K$
can be expressed in terms of the state vector $x$ via $K = K_0 + Mx$.

\begin{algorithm}
  \caption{Closed-Loop Transfer Function Tracking}
  \begin{algorithmic}[1]
  		\State \textbf{Output:} Optimal PID gains $K$
  		\State Hyperparameters: stepsize $\alpha > 0$, standard deviation $\sigma > 0$ for policy exploration, number of sampling directions $N$, vector of sampling times $T$
 		\State Initialize: PID gains vector $K_0$ of length $k$, policy $M = 0_{k \times n}$
 		\State Initialize: Target data $\tilde{x}$ for times in $T$
 		\State Set $K = K_0$
		\For{{each episode}}{}{}
        \State $x \leftarrow \texttt{step}(\texttt{feedback}(P,K),T)$
		\For {{each $j \text{ in } 1, \ldots, N$}}{}{}
		\State Choose $\delta_j \in \mathbb{R}^{k \times n}$ at random.
		\label{line:randomdelta}
		\State Perturb policy:
		\label{line:exploration}
        		\[
		        \begin{cases}
		            K^{+} \leftarrow (M + \sigma \delta_j)x + K_0\\
                    K^{-} \leftarrow (M - \sigma \delta_j)x + K_0
                 \end{cases}
                \]
        \State Calculate closed-loop step responses:
                \[
                \begin{cases}
                    x_{j}^{+} \leftarrow \texttt{step}(\texttt{feedback}(P,K^{+}),T)\\
                    x_{j}^{-} \leftarrow \texttt{step}(\texttt{feedback}(P,K^{-}),T)
                \end{cases}
                \]
        \State Evaluate rewards due to perturbation:
                \[
                \begin{cases}
                r_{j}^{+} \leftarrow \text{Reward at } x_{j}^{+}\\
                r_{j}^{-} \leftarrow \text{Reward at } x_{j}^{-}
                \end{cases}
                \]
        \EndFor
        \State Evaluate standard deviation $\sigma_r$ of the $2N$ rewards
        \State $M \leftarrow M + \frac{\alpha}{\sigma_r} \frac{1}{N} \sum_{j = 1}^{N} \big[ r_{j}^{+} - r_{j}^{-} \big] \delta_{j}$\label{eq:polUpdate}
        \State $K \leftarrow M x + K_0$\label{line:update}
        \EndFor
     \end{algorithmic}
    \label{alg:MainAlg}
\end{algorithm}

Intuitively, Algorithm \ref{alg:MainAlg} is exploring the parameter space of PID gains centered at $K_0$; note that $K_0$ remains fixed throughout training.
If we initialize $K_0$ with a zero vector then the first iteration simply operates with a zero-mean Gaussian matrix for the policy at line \ref{line:exploration}. 
Alternatively, if a set of PID gains is known to be stabilizing, or obtained through methods such as relay tuning or IMC, then we use those parameters to define $K_0$, and then evolve and improve them (guided by the reward objective) through evaluating the performance of small perturbations to the policy.

In line \ref{line:randomdelta}, the random matrices are determined by drawing each
entry from an independent standard normal distribution.
After the rewards are collected for several perturbed policies, in line \ref{eq:polUpdate} we update the policy using a scaled average of the finite-difference approximation for the gradient given by $\big[ r^{+} - r^{-} \big] \delta$. We scale by a fixed step-size $\alpha > 0$ as well as the reciprocal of the standard deviation of the rewards obtained from each sampled direction $\delta$. Together, these scaling factors give smaller update steps when the current policy is more sensitive to perturbations, while larger steps are permitted when the perturbed policies yield similar rewards. Finally, line \ref{line:update} updates $K$ on the basis of the new policy and current output measurements. The small number of hyperparameters illustrates the simplicity and interpretability of the
algorithm.

\subsection{Differences from other RL-based tuning approaches}

Here we highlight a few RL-based PID tuning approaches across various applications. With applications to wind turbine control, \cite{sedighizadeh2008adaptive} propose an actor-critic approach in which the PID gains are the actions taken by the actor at each time-step. \cite{carlucho2017incremental} develop an on-line discretization scheme of the state and action spaces, allowing for the implementation of the $Q$-learning algorithm for control of mobile robots. Finally, \cite{brujeni2010dynamic} implement the classical SARSA algorithm for control of chemical processes. Their approach uses IMC to define a collection of PID gains which comprises the action space. At each time-step in the control of a physical continuous stirred tank heater the algorithm selects the best gains.

Our approach does not require training a neural network for value function approximation nor to represent an actor. Instead, our policy is given by a matrix whose size is determined by the number of sampling times in an episode and the number of tunable parameters for a linear controller (e.g., PI or PID). Further, our policy update procedure occurs on a different time scale than the sampling time. In particular, we update the policy based on entire closed-loop step responses, rather than at each time-step of a step response. This distinction avoids an important phenomenon associated with switching control strategies. Namely, if two controllers are known to be stabilizing, switching between them can still destabilize the closed-loop (see Example 1 in \cite{malmborg1996stabilizing}). Closed-loop transfer function tracking is then an intuitive approach for embedding performance specifications into a reward function without destabilizing the closed-loop with stabilizing controllers. With this view, it is justified to treat PID parameters as actions in the RL framework.

\section{Simulation Results}\label{sec:results}

We present several simulation examples to illustrate our algorithm. The first example is a proof of concept in which we initialize the PID controller with unstable weights and construct a solution for Algorithm \ref{alg:MainAlg} to find. 
The second example initializes the PID parameters with SIMC \citep{skogestad2001probably}, then Algorithm \ref{alg:MainAlg} updates the PID parameters to compensate for slow changes in the plant gain. See Appendix \ref{app:details} for the hyperparameters used in Algorithm \ref{alg:MainAlg}.

\subsection{Example 1}

In this example, we demonstrate our tuning method by constructing a desired closed-loop transfer function with a given plant model and set of target PID parameters. 

Consider the following continuous-time transfer function:
\begin{equation}
P(s) = \frac{1}{(s+1)^3}.
\label{eq:ex1}
\end{equation}

\begin{figure}[thb]
\begin{center}
\includegraphics[width=8.4cm]{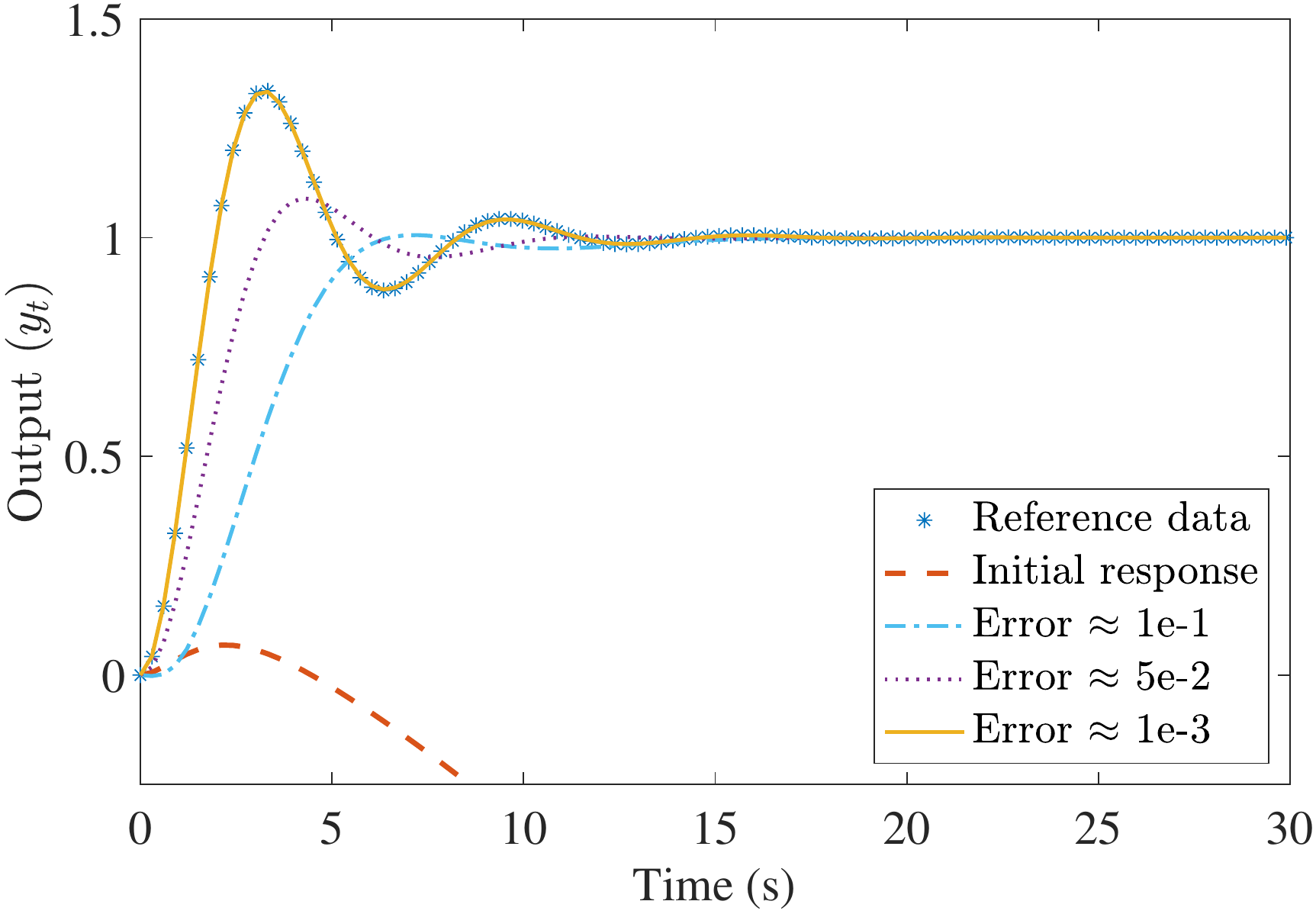}    
\caption{The closed-loop step response at the beginning is shown with a dashed line, at end of training shown with a solid line, along with the reference data.} 
\label{fig:ex1_fig1}
\end{center}
\end{figure}

We randomly initialize $k_p, k_i, k_d$ around zero and set the desired parameters to be $k_p = 2.5, k_i = 1.5, k_d = 1.0$. The initial parameters may destabilize the plant as shown in figure \ref{fig:ex1_fig1}. The target data then comes from uniform samples of the step response from the closed-loop transfer function $CP/(1 + CP)$ where $C$ is the PID controller with the aforementioned target parameters.

\begin{figure}[thb]
\begin{center}
\includegraphics[width=8.4cm]{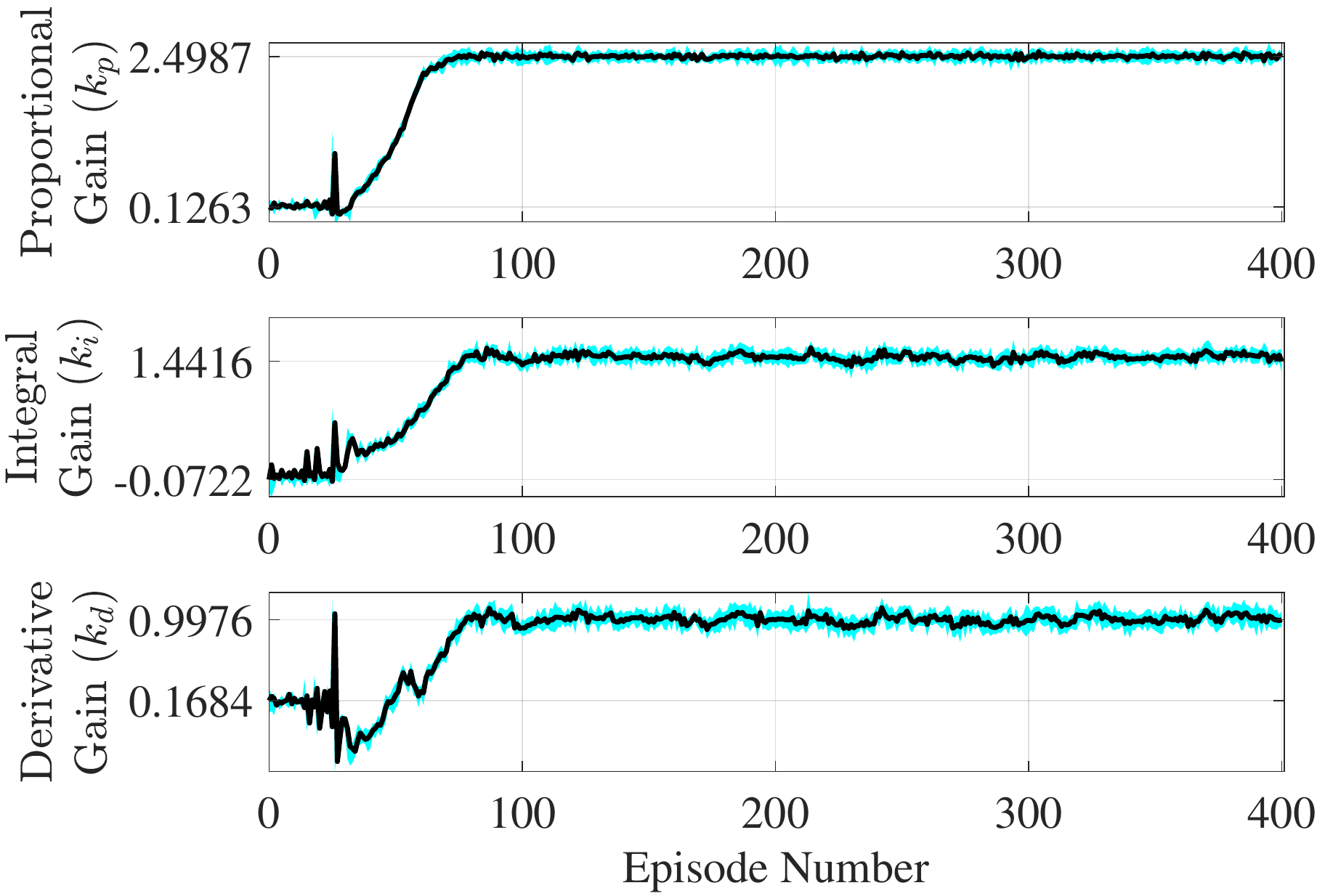}    
\caption{The value of the updated PID gains at each iteration is shown in black. The highlighted region shows the range of values seen at each episode in line \ref{line:exploration} of Algorithm \ref{alg:MainAlg}}
\label{fig:ex1_fig2}
\end{center}
\end{figure}

\begin{figure}[thb]
\begin{center}
\includegraphics[width=8.4cm]{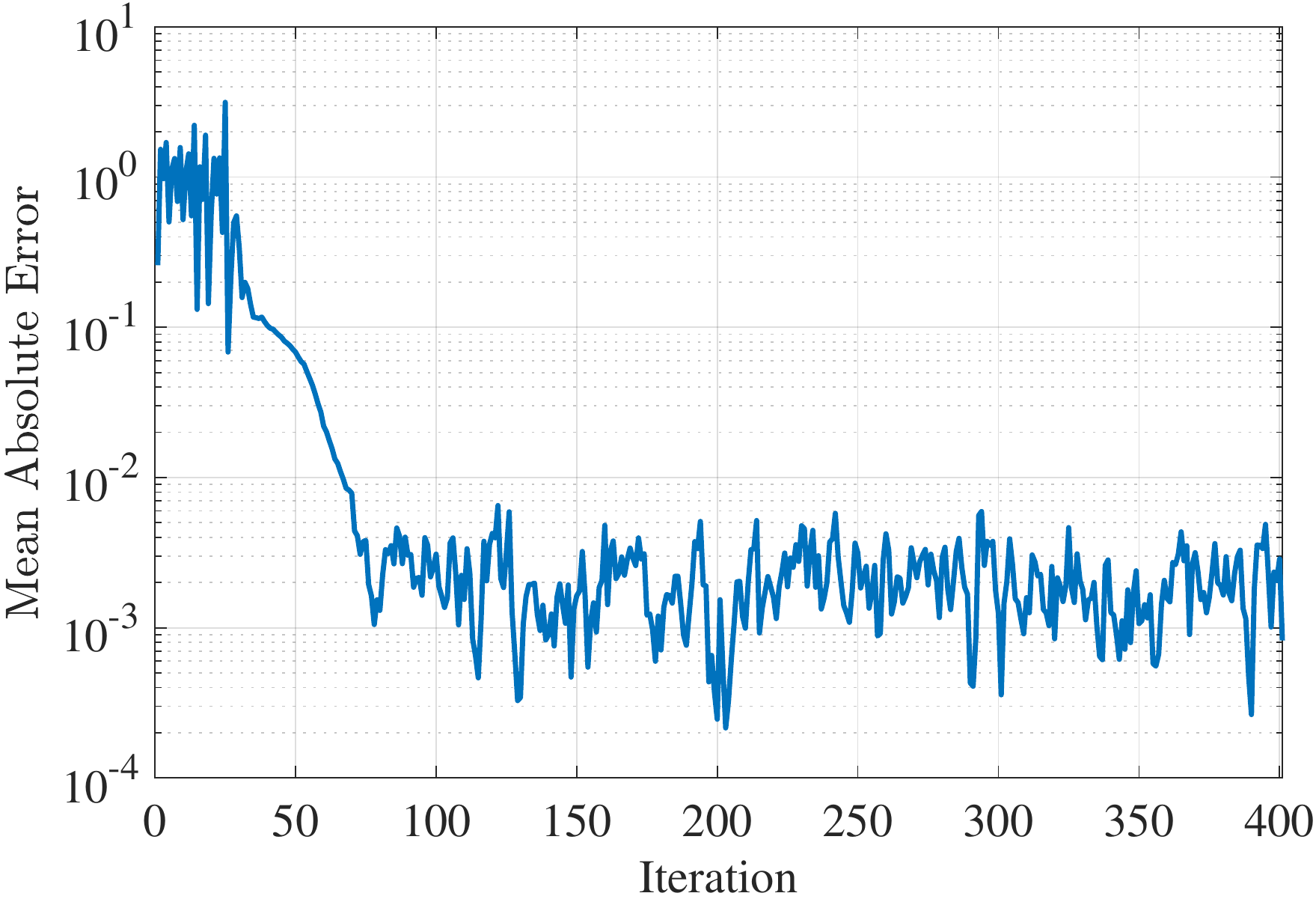}    
\caption{The mean absolute error decreases on a logarithmic scale with the number of episodes.} 
\label{fig:ex1_fig3}
\end{center}
\end{figure}

We highlight several important notes about this experiment. First, the speed at which the algorithm finds the correct parameters is determined by the step-size $\alpha$, the exploration parameter $\sigma$, and finally the relative distance between initial and target gains. We initialized the gains far away from the target to illustrate the trajectories of the gains during the learning procedure. Figure \ref{fig:ex1_fig2} shows the evolution of the PID parameters over the course of training. Note that the highlighted region indicates the parameters seen during the exploration described in line \ref{line:exploration} of Algorithm \ref{alg:MainAlg}. We run the simulation for many episodes to show the parameters hovering steadily around the solution. This behaviour aligns with the error (reward scaled by $-1$) curve  shown in figure \ref{fig:ex1_fig3}. We show several output responses in figure \ref{fig:ex1_fig1} corresponding to various levels of the error curve.

Our second remark is that the algorithm does not use any knowledge about the plant dynamics nor does it utilize a modeling procedure. Further, the PID control structure is only implicitly used, meaning the actions $K$ directly influence the closed-loop, but could correspond to any linear controller. Finally, the target step response is user-specified, which makes Algorithm \ref{alg:MainAlg} amenable to performance specifications. 

\subsection{Example 2}

In this example, we tune a PID controller using Algorithm \ref{alg:MainAlg} subject to drift in the plant gain and uncertainty in the time-delay.

Consider a nominal plant given by
\begin{equation}
P(s) = \frac{-0.02}{s + 1}e^{-s}.
\label{eq:ex4}
\end{equation}

We generate our target closed-loop data with $P$ under the PID gains obtained with the SIMC tuning method referenced in section \ref{sec:IMC}. These gains are $K_0$ in Algorithm \ref{alg:MainAlg}. We use SIMC for initialization due to its simplicity as well as to illustrate the compatibility of Algorithm \ref{alg:MainAlg} with existing tuning methods. Note that we are not comparing our algorithm against SIMC.

\begin{figure}
\begin{center}
\includegraphics[width=8.4cm]{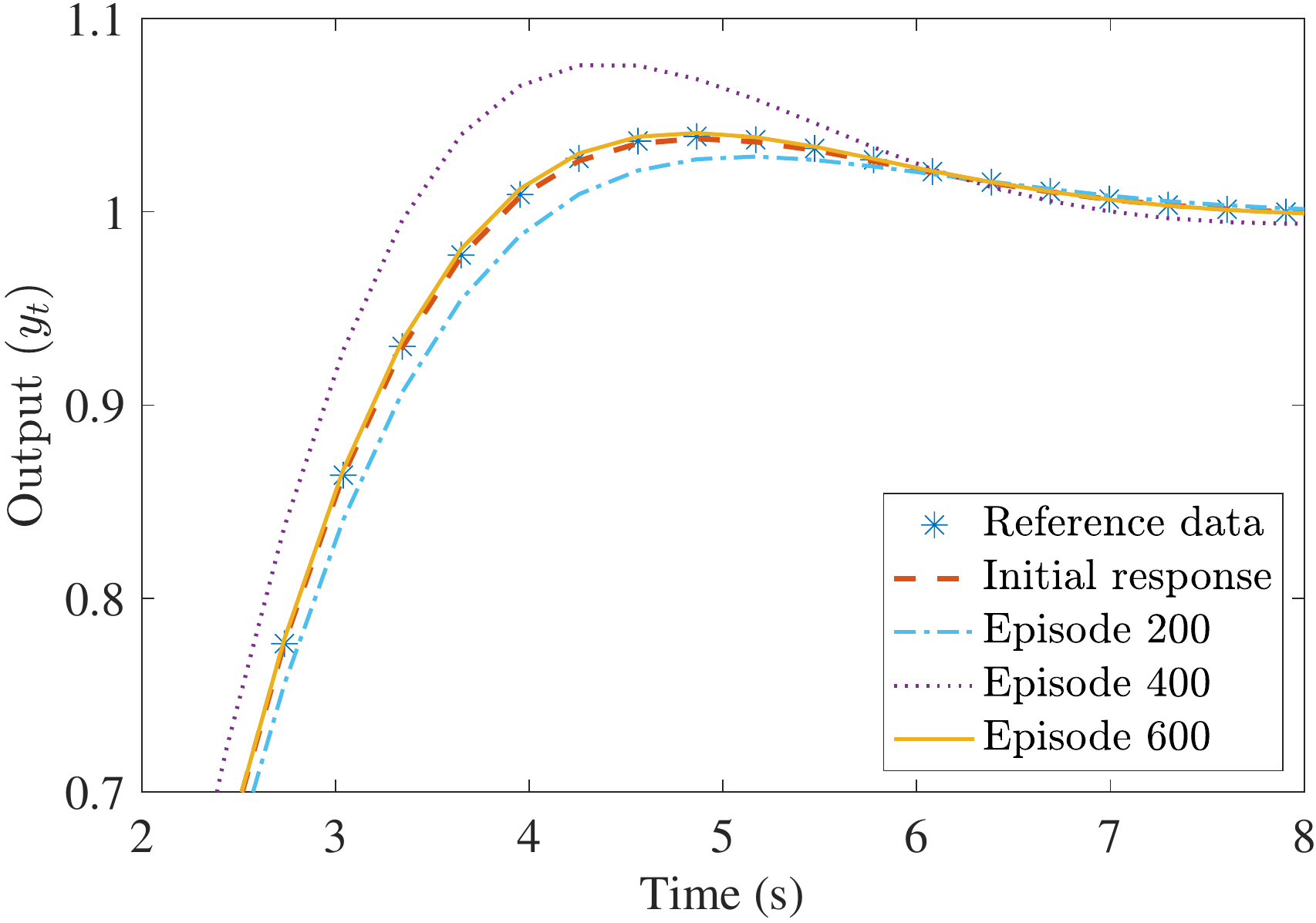}    
\caption{The closed-loop step response corresponding to different plant gains and adjusted PID parameters so as to maintain initial performance.} 
\label{fig:ex4_fig1}
\end{center}
\end{figure}

\begin{figure}
\begin{center}
\includegraphics[width=8.4cm]{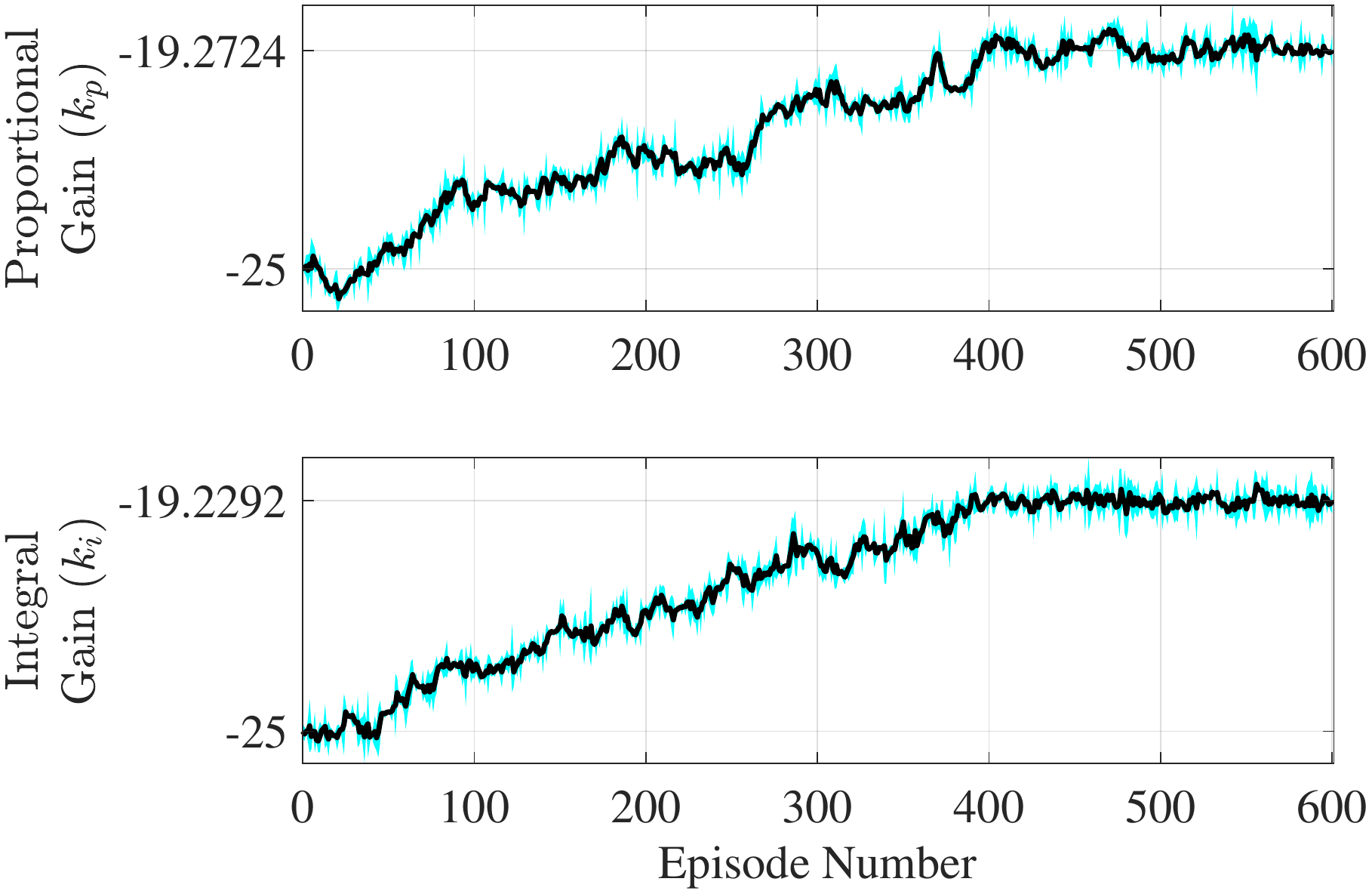}    
\caption{The updated integral and proportional gains at each episode. The plant gain remains fixed after episode 400.} 
\label{fig:ex4_fig2}
\end{center}
\end{figure}

\begin{figure}
\begin{center}
\includegraphics[width=8.4cm]{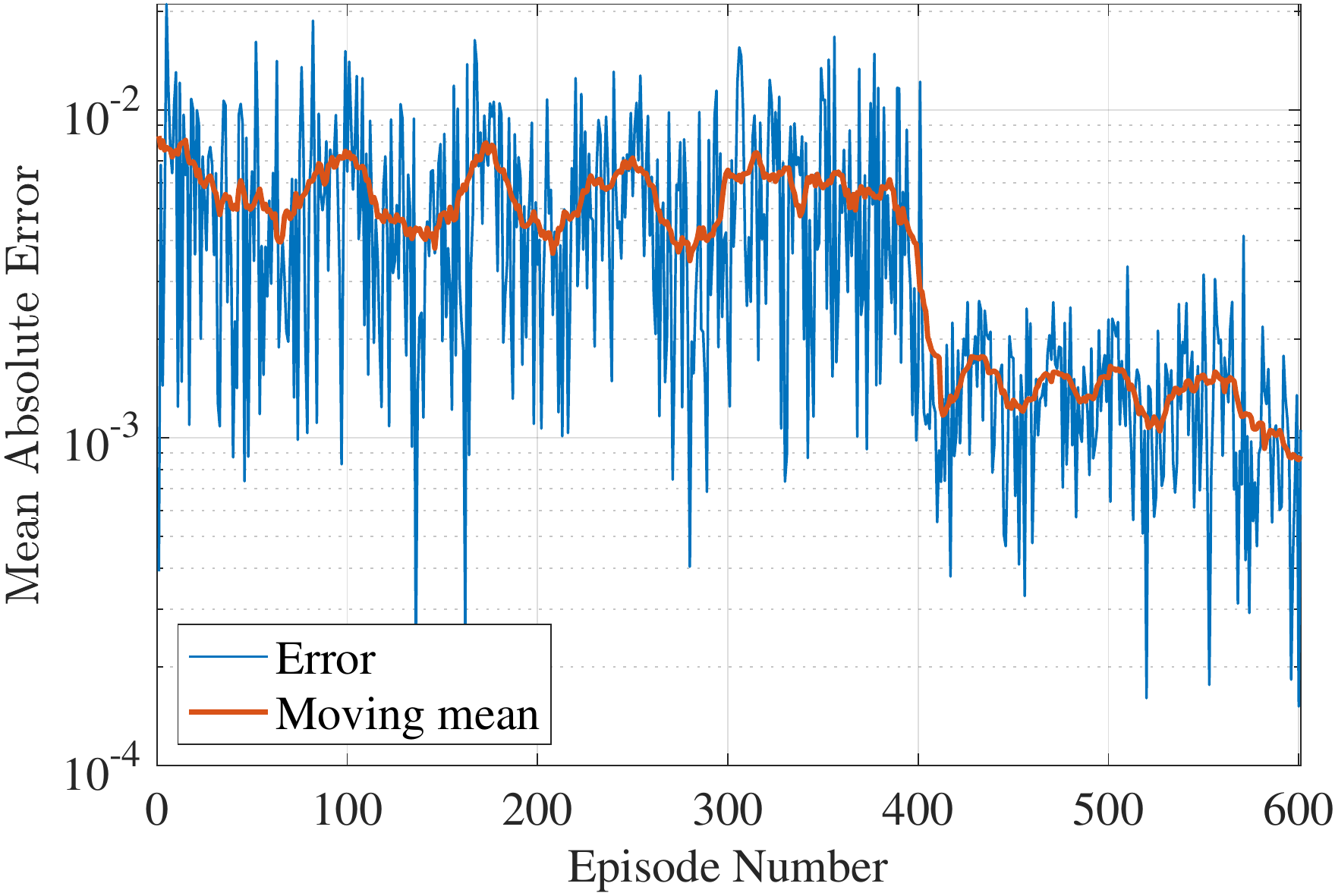}    
\caption{The mean absolute error at each episode is shown in blue with a red curve overlapping it to show the moving average across 20 episodes.} 
\label{fig:ex4_fig3}
\end{center}
\end{figure}

At the beginning of each episode we slightly change the gain in the numerator of $P$; we also perturb the time-delay from $1$ by adding a small amount of mean-zero Gaussian noise. The gain drifts linearly, so that at episode 400 its magnitude has increased by 30\%. (The final numerator is $-0.026$.) At episode 400 we keep the plant gain fixed simply to observe the parameter updates in figure \ref{fig:ex4_fig2} for both changing and static plant gains. Figure \ref{fig:ex4_fig3} shows the error being maintained as the plant gain drifts and figure \ref{fig:ex4_fig1} shows snapshots of the closed-loop output response throughout training. It is worth mentioning that the error is steadily being maintained even at the beginning of training. This is due to the SIMC-based initialization and local parameter improvements made by Algorithm \ref{alg:MainAlg}.



\section{Conclusion}

We have developed a simple and intuitive algorithm based on random search for tuning linear fixed-structure controllers. In the RL framework, we treat the entire closed-loop as the environment subject to new controller parameters as the actions. The reward function encodes performance requirements by considering a desired closed-loop transfer function. The simplicity of our algorithm allows for minimal hyper-parameter tuning, as well as straightforward initialization of the policy around a set of initial controller parameters.

\begin{ack}
We would like to thank Profs. Benjamin Recht and Francesco Borrelli of University of California, Berkeley for insightful and stimulating conversations. We would also like to acknowledge the financial support from Natural Sciences and Engineering Research Council of Canada (NSERC) and Honeywell Connected Plant.
\end{ack}
\small
\bibliography{ifacconf}             

\normalsize
\appendix
\section{Implementation Details}\label{app:details}    

We scripted Algorithm \ref{alg:MainAlg} in \texttt{MATLAB} and simulated the processes using the Control System Toolbox. 
We use different hyper-parameters for each example. However, we note that any set of hyper-parameters listed below lead to similar results for the each examples, but do not illustrate the parameter updates as clearly.\\
Example 1: $\alpha = 0.005$, $\sigma = 0.005$, $N = 10$\\
Example 2: $\alpha = 0.01$, $\sigma = 0.05$, $N = 1$ (2 policy perturbations per episode)\\
For all examples, samples were taken in increments of $0.30$ seconds.\\
It is also possible to incorporate momentum in the policy update (Line \ref{eq:polUpdate} of Algorithm \ref{alg:MainAlg}). This can lead to smaller and more steady errors and smoother parameter updates, but if the initial policy is unstable (e.g., Example 1) it can also exacerbate instability. We therefore omit it for simplicity.

\end{document}